\documentclass[12pt, smallcondensed]{scrartcl}
\usepackage{amsmath, amssymb, amsthm}
\usepackage{graphicx}
\usepackage{subfigure}
\usepackage{color}
\usepackage{cite}

\newcommand{\be}{\begin{equation}}
\newcommand{\ee}{\end{equation}}

\theoremstyle{plain}
\newtheorem{thm}{Theorem}[section]

\newtheorem{definition}{Definition}[section]

\title{Vertex-Transitive Polyhedra of Higher Genus, I}

\author{Undine Leopold}
              
\date{}

\begin{document}
\maketitle
\begin{abstract}
Since Gr\"unbaum and Shephard's investigation of self-intersection-free polyhedra with positive genus and vertex-transitive symmetry in 1984 the question of complete classification of such objects in Euclidean $3$-space has been open. 
Due to a recent article by G\'evay, Schulte, and Wills, we now know that the genus range $\mathfrak{g}\geq 2$ can only admit a finite number of vertex-transitive polyhedra, all with rotational Platonic symmetry. In this article, we show that the symmetry group must also act simply transitively on the vertices. Furthermore, the case of rotational tetrahedral symmetry is settled completely - the unique example is already known.\\[1ex]
\end{abstract}
\section{Motivation}\label{sec:motivation}

Vertex-transitive polyhedra of positive genus are generalizations of the well-known convex uniform polyhedra. Examples of non-spherical vertex-transitive polyhedra were first examined by Gr\"unbaum and Shephard in a $1984$ survey article on \emph{Polyhedra with transitivity properties} \cite{gs1984}. Besides two two-parameter families of toroidal vertex-transitive polyhedra, only five examples of genus $\mathfrak{g}\geq 2$ were given. To this day it is unclear whether all such polyhedra of genus two and higher have been found. A sixth example, the combinatorially regular Gr\"unbaum polyhedron, has been discovered several times. Recent progress has been made by G\'evay, Schulte, and Wills \cite{gsw2013}, not only in providing the seventh example, but in restricting the symmetry groups to the rotation groups of the Platonic solids. Yet the principal question of complete enumeration of the combinatorial types is still open. 

Throughout this article, a \emph{polyhedron} $P$ is a closed, connected, orientable surface embedded in $\mathbb{E}^3$, tiled by finitely many plane simple polygons (\emph{faces}) in a face-to-face manner. The polygons are called \emph{faces} of $P$, their vertices and edges are called the \emph{vertices} and \emph{edges} of $P$, respectively. Gr\"unbaum \cite{gruenbaum1999} called such polyhedra \emph{acoptic}. A polyhedron $P$ can be considered a \emph{polyhedral realization}, that is, an embedding into $3$-space, of an underlying \emph{polyhedral map} $M$ on a closed orientable surface $S$. We refer the reader to the survey articles \cite{bs1997} and \cite{bw1993} for background on polyhedral maps and polyhedral manifolds.

Associated with a polyhedron $P$, we have a combinatorial symmetry group $\Gamma$ (automorphisms of $M$) and a geometric symmetry group $G$ (linear isometries of $\mathbb{E}^3$ fixing $P$ and inducing an automorphism in $\Gamma$). 
$P$ is \emph{(geometrically) vertex-transitive} if $G$ acts transitively on the vertices of $P$. Since these polyhedra are free of non-trivial self-intersection, the problem at hand is distinct from the classification of the $75$ uniform polyhedra by Coxeter, Longuet-Higgins, Miller \cite{coxeter1954}, Sopov \cite{sopov1970}, and Skilling \cite{skilling1975}, and the enumeration of regular, uniform, or equivelar maps on surfaces \cite{babai1991,thomassen1991,bk2008,conder2009,pellicer2011,kn2012,pellicer2014} with its abundance of examples. In fact, we investigate in the spirit of \cite{mcm1983,sw1985,wills1986,brehm1987,bokowski1989,sw1991,ssw2002,sw2012}, requiring a particular high geometric symmetry yet not necessarily combinatorial regularity. 

Our goal is to establish the results in Theorem \ref{thm:simple} and Theorem \ref{thm:tet}.
\begin{thm}\label{thm:simple}
Let $P$ be a polyhedron of genus $\mathfrak{g}\geq 2$ and $G$ its geometric symmetry group, such that $G$ acts transitively on the vertices of $P$. Then $G$ must act simply transitively on the vertices of $P$.
\end{thm}

\begin{thm}\label{thm:tet}
There is exactly one combinatorial type of polyhedron of genus $\mathfrak{g}\geq 2$ which can be embedded into $3$-space with vertex-transitive tetrahedral symmetry. This polyhedron has genus $3$ and is one of the known examples in \cite{gs1984} (based on a snub tetrahedron, and with a combinatorial snub tetrahedron as its convex hull). The polyhedron has $12$ vertices, $48$ edges and $32$ triangular faces and is of Schl\"afli type $\{3,8\}$, which means that $8$ triangles meet at each vertex.
\end{thm}

The article is organized as follows. After proving Theorem \ref{thm:simple} in Section \ref{sec:simple}, we focus on the case of tetrahedral rotation symmetry $T$ in the remaining sections. 
We enumerate possible maps with an orientation-preserving simply vertex-transitive action of $T$, up to \emph{geometric isomorphism}. 
The list of maps which we actually have to check for realizability turns out to be rather short. Finally, only one polyhedron of tetrahedral symmetry emerges as realizable - it is of genus $3$ and part of Gr\"unbaum and Shephard's \cite{gs1984} list. We provide coordinates for its vertices, which appears to not have been done previously. Furthermore, since these coordinates must satisfy certain polynomial inequalities, we are able to determine them among the small integers. We give combinatorial data as well as integer coordinates for the vertices in Table \ref{tab:map1coord}.  The polyhedron is depicted in Figures \ref{fig:tet_poly} and \ref{fig:tet_poly2} at the end of the article. 

\section{Preliminaries}\label{sec:intro1}

The well-known polyhedron formula of Euler relates the two invariants genus $\mathfrak{g}$ and Euler characteristic $\chi$ of the map $M$ underlying a polyhedron $P$ to the combinatorial data:
\be\label{euler}\chi=2-2\mathfrak{g}=|V|-|E|+|F|,\ee 
where $|V|$ is the number of vertices, $|E|$ is the number of edges, and $|F|$ is the number of faces. In a \emph{triangulated} polyhedron of genus $\mathfrak{g}$ with $|V|$ vertices, that is a polyhedron with only triangular faces, we have $|F|=\frac{2}{3}|E|$. 
The same holds for any (polyhedral) map consisting only of triangular faces, a so-called \emph{triangulation} of the surface.
Therefore, we obtain for the number of edges $|E|$ in a \emph{triangulated} surface of genus $\mathfrak{g}$
\be\label{minedges}|E|=3|V|+6(\mathfrak{g}-1).\ee
A triangulation has the highest number of edges among all (polyhedral) maps with a given number of vertices on the surface of genus $\mathfrak{g}$. 
Note that if we bound this maximal number of edges above by $|V|(|V|-1)/2$ we obtain the Heawood inequality \cite{heawood} which, more generally, provides a lower bound for the number of vertices of a polyhedral map on a surface of genus $\mathfrak{g}$
\be \label{heawood}|V|\geq \left\lceil \frac{1}{2}\left(7+\sqrt{49+48(\mathfrak{g}-1)}\right)\right\rceil=\left\lceil \frac{1}{2}\left(7+\sqrt{1+48\mathfrak{g}}\right)\right\rceil.\ee
This inequality is sharp for all orientable surfaces except the surface with genus $\mathfrak{g}=2$, see \cite{ringel1980}.

 All combinatorial types of vertex-transitive polyhedra of genus zero are known; they comprise the Platonic solids, the Archimedean solids, and the infinite families of prisms and antiprisms. Remarkably, two two-parameter families exist for genus one (see \cite{gs1984} and \cite{gsw2013} for descriptions). In this article, we focus our attention on $\mathfrak{g}\geq 2$ as only a few examples of higher genus are known, and the completeness of the list has never been established. Moreover, G\'evay, Schulte, and Wills recently proved the following two theorems.

\begin{thm}\label{thm:sym}
(\cite[Theorem 4.1]{gsw2013}) There are only finitely many vertex-transitive polyhedra in the genus range $\mathfrak{g}\geq 2$. The symmetry group of each vertex-transitive polyhedron in the genus range $\mathfrak{g}\geq 2$ is a Platonic rotation group.
\end{thm}

In particular, $G$ cannot contain reflections. Recall that the Platonic rotation groups are $T$ (tetrahedral rotation group), $O$ (octahedral rotation group), and $I$ (icosahedral rotation group). 
 We know further that all vertices of $P$ lie on a sphere (around the point fixed by the finite symmetry group), and that all faces of $P$ are convex \cite{gs1984}.

\begin{thm}\label{thm:genus01}
(\cite[Theorem 4.2]{gsw2013}) A vertex-transitive polyhedron with a reducible symmetry group must have genus $0$ or $1$. There are infinitely many vertex-transitive polyhedra of genus $0$, as well as of genus $1$.
\end{thm}

Tucker \cite{tucker2013} investigated how the genus of an embedded smooth surface limits its possible geometric symmetries. 
The most prominent limitation is that an axis of rotational symmetry must pierce an embedded surface in an even number of points (possibly none). This and some other observations gave rise to the following theorem in \cite{tucker2013}, with a minor correction in the joint article \cite{lt2014}.

\begin{thm}\label{thm:genus} (Tucker \cite{tucker2013}, \cite{lt2014}) The following are necessary and sufficient conditions on the genus $\mathfrak{g}$ of a surface in order to admit a smooth embedding into $\mathbb{E}^3$ with the listed symmetry. 
\begin{itemize}
\item Rotational tetrahedral symmetry is possible if and only if \[\mathfrak{g}=m_1\cdot6+m_2\cdot8+\{0,3,5,7\},\]
\item Rotational octahedral symmetry is possible if and only if \[\mathfrak{g}=m_1\cdot 12+m_2\cdot16+m_3\cdot18+\{0,5,7,11,13\},\]
\item Rotational icosahedral symmetry is possible if and only if \[\mathfrak{g}=m_1\cdot 30+m_2\cdot 40+m_3\cdot 48+\{0,11,19,21,29,31,37\},\]
\end{itemize}
for some nonnegative integers $m_1$, $m_2$, $m_3$.

\end{thm}

We conclude from this theorem that the genera of our polyhedra must be compatible with the chosen symmetry, even though polyhedra are not smooth. Smoothness was assumed in order to make a more general statement, avoiding pathological embeddings, but the statements remain true for surfaces embedded as polyhedra (which can easily be smoothed out while keeping the vertex coordinates and remaining on an embedded surface). For the purposes of this article, it is not necessary to do a detailed study of the Riemann--Hurwitz equation as in \cite{lt2014}.

\subsection{Maximal Triangulation}

When a face of a vertex-transitive polyhedron $P$ under a Platonic rotation group $G$ has a trivial stabilizer, this face can be triangulated by some of its diagonals in a non-intersecting way, and the action of the symmetry group takes the triangulation of this face to the other faces in the same orbit. In this way an orbit of $|G|$ $m$-gons (with trivial stabilizer) is split into $m-2$ orbits of $|G|$ triangles each. Observe that no new vertices are added.

If a face of $P$ has a non-trivial stabilizer under $G$, this face is fixed by a rotation of some order $n$, which is the generator of the stabilizer subgroup, such that the face is a convex $k\cdot n$-gon (for some positive integer $k$, and $n=2,3,4,5$ depending on the symmetry group, see Theorem \ref{thm:sym}). Now it is possible to split this face by diagonals, except for when $n=2$, into a ``central'' regular $n$-gon and $(k-1)n$ adjoining triangles which lie symmetrically with respect to the axis of the stabilizer. For $n=2$ (and necessarily $k>1$) it is possible to triangulate with a diagonal meeting the axis of the stabilizing rotation, and other diagonals placed symmetrically with respect to the axis. Again, the symmetry group of the polyhedron then takes care of the splitting of the other faces in the orbit of the original face. Thus, the orbit of $\frac{|G|}{n}$ $k\cdot n$-gons gives rise to one orbit of $\frac{|G|}{n}$ $n$-gons (if $n>2$) and $k-1$ orbits consisting of $|G|$ triangles. 

As a result of the above modifications of the face structure of $P$, we obtain a vertex-transitive polyhedron $P'$ whose only non-triangular faces may be squares (in the case of octahedral symmetry) or regular pentagons (in the case of icosahedral symmetry). We call a polyhedron $P$ with such a face structure a \emph{maximally triangulated} polyhedron, even though some of its faces may not be triangles.

 Our investigation of vertex-transitive polyhedra, as well as the proof of Theorem \ref{thm:simple}, is greatly simplified by assuming $P$ to be a maximally triangulated polyhedron. When no maximally triangulated polyhedron can be found, then the search for less triangulated versions obtained by merging faces is futile. However, if a maximally triangulated polyhedron is found, it is still possible to investigate whether certain faces can be made coplanar and combined into larger faces while preserving vertex-transitivity.

\section{Proof of Theorem \ref{thm:simple}}\label{sec:simple}
The proof of simple transitivity of the symmetry group on the vertices of a vertex-transitive polyhedron is by cases. 
\subsection{Tetrahedral Symmetry}\label{sec:tet:simple}
We first prove simple transitivity for the case that $G=T$, the rotation group of a regular tetrahedron. Since $T$ acts transitively on the vertices, there are necessarily $12$, $6$, or $4$ vertices, depending on whether a vertex-stabilizer is trivial, or generated by a rotation of order $2$ or $3$. 
Polyhedra on surfaces of genus $\mathfrak{g}\geq 2$ have at least $10$ vertices by the Heawood inequality \eqref{heawood}, necessitating that 
\[|V|\geq \left\lceil \frac{1}{2}\left(7+\sqrt{1+48\mathfrak{g}}\right)\right\rceil\geq\left\lceil \frac{1}{2}\left(7+\sqrt{97}\right)\right\rceil=9,\]
as well as the works of Ringel and Jungerman \cite{ringel1980} who show that $9$ do not suffice.
We conclude that there are $12$ vertices in the single vertex orbit, each of which is trivially stabilized under $T$. This means that $T$ acts simply transitive on any higher genus vertex-transitive polyhedron with tetrahedral symmetry.

\begin{figure}[h] 
\begin{center}
\includegraphics[width=.45\textwidth]{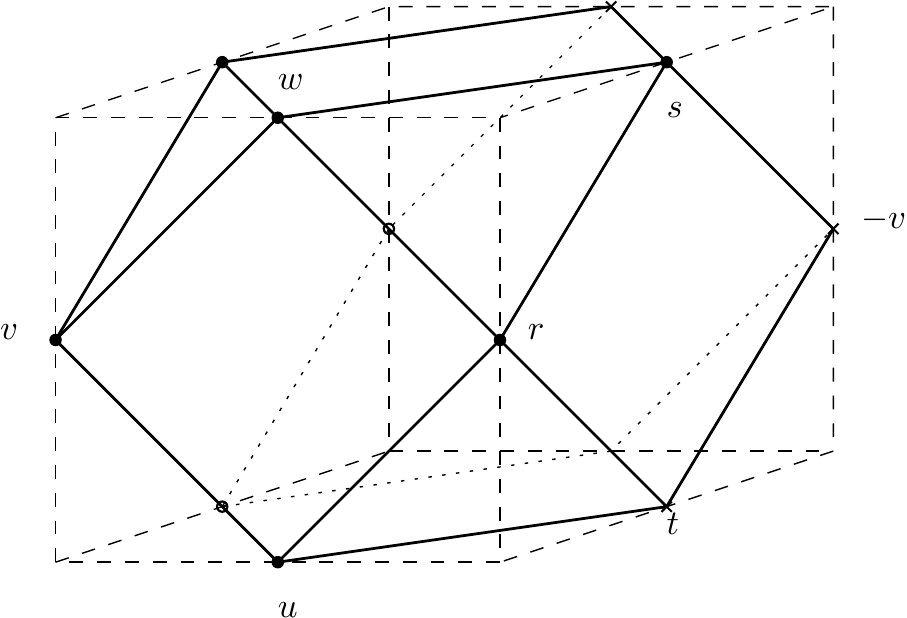}
\end{center}
\caption{Ruling out a cuboctahedron as convex hull.}
\label{fig:plat1}
\end{figure}

In addition, consider the auxiliary cube as in Figure~\ref{fig:plat1}, such that the order $2$ axes of the tetrahedral group pierce its face centers and the order $3$ axes are its space diagonals. If the vertices of $P$ are the midpoints of the edges of the auxiliary cube, i.e. $P$ has the cuboctahedron as convex hull, an initial vertex $v$ could connect to all vertices which are neighbors on this convex hull, such as $w$ and $u$, giving two edge orbits of size $12$ and a total of $12+12=24$ edges. Furthermore, $v$ could be joined by an edge to precisely one of the other vertices on the incident square faces of the convex hull, such as $r$, but not to its opposite $-v$, without producing non-trivially intersecting edges (all edges in the same orbit would have to go through the origin), giving another $6$ edges. Of the remaining two theoretically possible edge orbits of size $12$, one represented by the segment $vs$, one represented by the segment $vt$, only one could occur, otherwise an intersection is produced by the diagonals of the trapezoid
$vws(-v)$. This means $P$ can have at most $12+12+6+12=42$ edges. Furthermore, if $P$ has tetrahedral symmetry and is not spherical ($\mathfrak{g}>0$), it has at least genus $3$ by Theorem \ref{thm:genus}. As a maximally triangulated surface, it would have to be completely triangulated. A triangulation on a surface of genus $\mathfrak{g}\geq 3$ has precisely $3|V|+6(\mathfrak{g}-1)$ edges (see Equation (\ref{minedges})), so we would need at least $3\cdot12+6\cdot(3-1)=48$ edges. Hence, the vertices of $P$ cannot lie on the edge midpoints of the auxiliary cube, and the convex hull of $P$ cannot be a cuboctahedron.

\subsection{Octahedral Symmetry}\label{sec:oct:simple}
For a vertex-transitive polyhedron $P$ with symmetry group $O$, the rotation group of a regular octahedron or cube such as the auxiliary cube in Figure~\ref{fig:plat1}, the vertex orbit can have size $6$, $8$, $12$, or $24$, depending on the size of the vertex-stabilizer ($4$, $3$, $2$, or $1$). As before in the tetrahedral case, fewer than $12$ vertices cannot support a polyhedron on a surface of higher genus $\mathfrak{g}\geq 2$. Twelve vertices can only occur if the vertices lie on the axes of order $2$. In this case the tetrahedral rotation group $T$, being a (geometric!) subgroup of the octahedral group $O$, still acts on the polyhedron and happens to act transitively on the vertices. However, in Section \ref{sec:tet:simple} we have established that these positions for the vertices are infeasible for a vertex-transitive polyhedron with tetrahedral symmetry. 
Consequently, any vertex-transitive polyhedron with octahedral symmetry and genus $\mathfrak{g}\geq 2$ must have precisely $24$ (trivially stabilized) vertices. Then, $O$ must act simply transitive on the vertices of $P$, which proves Theorem \ref{thm:simple} for the case of octahedral symmetry.

\subsection{Icosahedral Symmetry}\label{sec:ico:simple}
In the case of rotational icosahedral symmetry $I$ we assume again that $P$ is maximally triangulated. Maximal triangulation yields either a completely triangulated polyhedron, or a polyhedron with some residual regular pentagon faces (pierced and stabilized by $5$-fold axes). Furthermore, a vertex-transitive polyhedron $P$ with icosahedral symmetry of genus $\mathfrak{g}\geq 2$ must actually have at least genus $\mathfrak{g}=11$ according to Tucker's result in Theorem \ref{thm:genus}. In theory, it could have $12$, $20$, $30$, or $60$ vertices, each with an order $5$, order $3$, order $2$, or trivial stabilizer, respectively. In practice, only $60$ vertices can occur. This can be seen as follows.

First, a vertex-stabilizer of order $5$ can be immediately excluded. It leads to a polyhedron with only $12$ vertices, too few for a polyhedron with genus $\mathfrak{g}=11$ or higher; the Heawood inequality (\ref{heawood}) stipulates $|V|\geq 15$ for $\mathfrak{g}\geq11$.

\begin{figure}[ht] 
\begin{center}
\includegraphics[width=.3\textwidth]{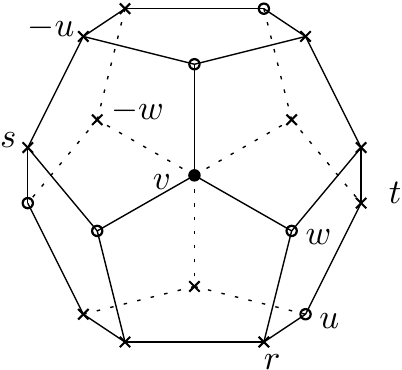}
\end{center}
\caption{Ruling out a dodecahedron as convex hull.}
\label{fig:plat2}
\end{figure}

Second, we show that $20$ vertices are not possible. Suppose $P$ has $20$ vertices. Then each vertex is stabilized by a rotation of order $3$, and ${\rm conv}(P)$ is a regular pentagonal dodecahedron. 
 $P$ can have at most $30+60=90$ edges if it is maximally triangulated, which we can easily see from the sketch of the convex hull in Figure~\ref{fig:plat2}.
 
Without producing self-intersections, an initial vertex $v$ could connect to its neighbors on the convex hull, such as $w$, giving exactly one edge orbit of size $30$ (each of these edges is fixed by a rotation of order $2$). The vertex $v$ cannot connect (via an edge) to its opposite $-v$ because edges cannot pass through the origin. Furthermore, $v$ cannot connect to the non-adjacent vertices in the incident pentagons, such as $s$, as the $5$-fold rotations which fix the pentagons produce self-intersections in that edge orbit. For the same reason, $v$ cannot connect to $-w$ and symmetrically placed vertices with respect to the axis through $v$ and $-v$ (for example, the edge between $v$ and $-w$ intersects with the symmetric edge between $r$ and $-u$). Finally, $v$ can only connect to one of $u$, $t$ (and symmetrically placed vertices) without producing self-intersections; if both corresponding edge orbits are in $P$, then so are the diagonals of the trapezoid $vwus$, and we have a self-intersection. The resulting edge orbit of either $vu$ or $vt$ has size $60$, giving a total of at most $90$ edges.

In order to triangulate $P$ completely, without preserving symmetry, we need at least \[3|V|+6\cdot(11-1)=60+60=120\] edges because of (\ref{minedges}) and the fact that $P$ has at least genus $11$ (see Theorem \ref{thm:genus} and recall that we assume $\mathfrak{g}\geq 2$). Since a maximally triangulated $P$ (i.e., triangulating as far as possible while keeping the vertex-transitive symmetry) has at most $90$ edges, it must have residual pentagonal faces. However, the only possible set of pentagonal faces, considering the possible edges incident to $v$ and recalling that these pentagons need a non-trivial stabilizer, are those which close up to give the convex hull, a dodecahedron of genus $0$. 
We conclude that  $20$ vertices are not sufficient for an (embedded) vertex-transitive polyhedron of higher genus with icosahedral symmetry.

\begin{figure}[h]
\begin{center}
\includegraphics[width=.3\textwidth]{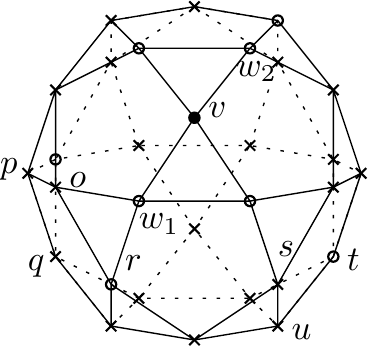}
\end{center}
\caption{Ruling out an icosidodecahedron as convex hull.}
\label{fig:plat3}
\end{figure}

Next we also reject the possibility that there are $30$ vertices. Suppose $P$ has $30$ vertices. Then each vertex is stabilized by a rotation of order $2$, and ${\rm conv}(P)$ is a fully vertex-truncated regular dodecahedron (icosidodecahedron), as in Figure~\ref{fig:plat3}. This solid's symmetry implies that a maximally triangulated $P$ can have at most $60+30+30=120$ edges. Again, we can explain this by examining the symmetries of the convex hull and deducing the intersections of edges.

As before, an initial vertex $v$ could connect to all its neighbors on the convex hull, such as $w_1$, giving one edge orbit of size $60$. Again, $v$ cannot connect to its opposite $-v$ as well as non-adjacent vertices, such as $o$, on the incident pentagons of the convex hull. Also, $v$ cannot connect to any other vertex, such as $p$, with the resulting edge meeting an order $5$ axis. The vertex $v$ could only connect to one of $r$, $s$; if both resulting edge orbits (of size $30$, since each edge is stabilized by a two-fold rotation) were allowed, we would obtain self-intersections, for example, from the diagonals of $vw_1rw_2$. Now observe that there are planar decagons (equatorial decacons of the icosidodecahedron) made up of edges of ${\rm conv}(P)$ going through $v$ and $r$, and through $v$ and $s$, respectively. It is clear that $v$ cannot connect to any of the remaining vertices on these decagons because the order $5$ rotations which fix the decagons would produce self-intersections. Furthermore, only one of the edge orbits containing $vq$ or $vt$ could be permitted, again because we want to avoid crossing diagonals in planar quadrilaterals such as the rectangle $vw_2tu$. This would give one more edge orbit of size $30$ (edges fixed by two-fold rotations), and completes our counting argument.

Triangulating $P$ entirely (without regard for symmetry), for $\mathfrak{g}\geq 11$, we need at least  \[3|V|+6\cdot10=90+60=150\] edges. However, we have just shown that a symmetry-compatible maximal triangulation has no more than $120$ edges. Thus a maximally triangulated $P$ must have faces of $5$-fold rotational symmetry pierced by $5$-fold axes. There are two kinds of non-trivially stabilized pentagon faces which are permitted by the edges incident to $v$. The first kind, regular pentagon faces given by alternate vertices of the equatorial decagons, will contain the origin, and therefore intersect in their interiors with their images under a three-fold rotation, say. The second kind of symmetric pentagon faces are the twelve pentagon faces of ${\rm conv}(P)$. However, using them implies that the completed (but now unsymmetrical) triangulation of $P$ (obtained by adding two edges in each pentagon) can have no more than $120+12\cdot2=144$ edges, a contradiction. In consequence, $30$ vertices are also not sufficient to construct a vertex-transitive polyhedron.

As with the tetrahedral and octahedral symmetry groups this only leaves the possibility of a trivial vertex stabilizer and a maximal number of vertices, here $60$. It is thus established that if $I$ acts vertex-transitively on a polyhedron $P$ of genus $\mathfrak{g}\geq 2$, then it must act simply transitively, and we have completed the proof of Theorem \ref{thm:simple} for all cases.

\section{Further Background}\label{sec:background2}

For the following treatment, we encode the geometric symmetry of a polyhedron directly into its underlying map which leads into the territory of topological graph theory. See \cite{topograph} for general background. 

Labeling each \emph{dart} (or \emph{arc}, directed edge) connecting two vertices $v$ (\emph{tail}), $w$ (\emph{head}) of a maximally triangulated vertex-transitive polyhedron $P$ under the group $G$ with the unique element $g \in G$ such that $(P)g=P$ and $(v)g=w$ leads to a labeling of the darts of the underlying polyhedral map. The opposite dart is necessarily labeled with $g^{-1}$, and so we may choose to display only one of the darts as a directed edge. 

Given a polyhedron $P$ and a choice of orientation, those darts forming an oriented (counterclockwise) boundary walk of a face $F$ are called \emph{associated} with $F$. The cyclic sequence of outgoing darts at any vertex $v$ is called its \emph{local rotation}, a term which we will also use for the cyclic sequence of associated dart labels. Observe that those labels do not repeat in the local rotation. By the process of labeling, we obtain the \emph{underlying labeled map}, which can be encoded succinctly (in either orientation) either by the local rotation at a chosen \emph{initial vertex} $v$ or by a list of so-called \emph{orbit symbols} with respect to $v$. Note that for the remainder of the article, the choice of orientation is an integral element of any (underlying) labeled map, so a switch in orientation produces a map distinct from the first; even though no dart labels have changed, the face associations of darts have changed. 

\begin{figure}
\begin{center}
\subfigure[Orbit of type 2.]{\includegraphics[width=.25\textwidth]{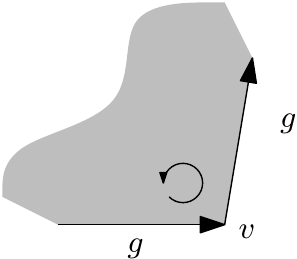}\label{fig:vertex:1}}\qquad
\subfigure[Orbit of type 1.]{\includegraphics[width=.5\textwidth]{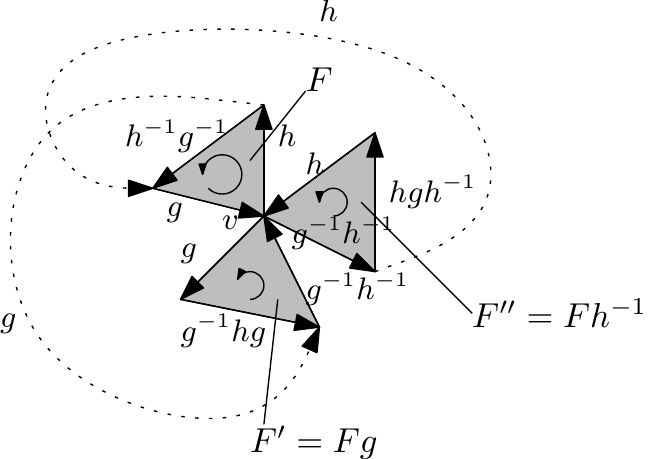}\label{fig:vertex:2}}\\
\end{center}
\caption{Faces of the same orbit at $v$.}
\label{fig:vertex}
\end{figure}

We note that all dart labels must be so-called \emph{core rotations} of $G$, that is, rotations around the corresponding axis with the smallest possible (positive or negative) non-trivial rotation angle. The set of core rotations is denoted $R(G)$. For example, half-turns around an order $4$ axis are not permitted; the corresponding edge necessarily intersects with an edge in the same orbit on the rotation axis. Next, the neighborhood of a chosen initial vertex $v$ partitions into faces from face orbits of the following two types as in Figure \ref{fig:vertex}. \emph{Type 1} is an orbit of trivially stabilized faces, whereas \emph{type 2} is an orbit of non-trivially stabilized faces. The former ones must be triangles since we assume maximal triangulation, whereas each of the latter faces have the property that their ${\rm ord}(g)>2$ associated darts bear the same core rotation $g$ as label (those faces in the polyhedron are then regular ${\rm ord} (g)$-gons).

The orbit symbol for an orbit of type 2, as in Figure \ref{fig:vertex:1}, is given by $(g)$ for the core rotation $g$ with ${\rm ord}(g)>2$ labeling the darts of the only face of this orbit incident at $v$. We are able to recover from it the \emph{oriented combinatorial angle} $(g,g)$ of this face at $v$, that is, the labels of consecutive darts associated to the face and joining at $v$.

For an orbit of type 1, there are three faces in the orbit incident at $v$, and three consecutive pairs of labels, or oriented combinatorial angles, can be found. If one is $(g,h)$, then the others are $(h,g^{-1}h^{-1})$ and $(g^{-1}h^{-1},g)$ by symmetry. Each of these pairs can be found as consecutive entries in the orbit symbol $(g,h,g^{-1}h^{-1})$. Recall that $g, h, g^{-1}h^{-1}$ are core rotations, and note that the orbit symbol is a cyclic symbol, i.e. $(h,g^{-1}h^{-1},g)$and $(g^{-1}h^{-1},g,h)$ encode the same orbit in the same orientation. A change of assigned orientation for the surface or map therefore inverts all entries of and reverses the cyclic order of entries in the orbit symbol; in this case we would obtain $(h^{-1},g^{-1},hg)$. 

It is worth mentioning that the underlying labeled map (a map on an oriented surface!) can be obtained as a branched cover from a labeled, oriented, one-vertex quotient map, such that the orbit symbols give precisely the labels in the oriented face boundary walks of the quotient. The term \emph{quotient map}, however, necessitates a more general, non-polyhedral notion of map: a 2-cell decomposition of a closed surface with precisely one vertex, allowing semi-edges. Semi-edges have only one vertex, and only one associated dart, and are the image of those edges upstairs with an involutory label. See also \cite{topograph,ga1974,gross1974,ga1973,jones1978}. 

\subsection{Candidate Maps}

With the encoding into face orbits, we call a collection of orbit symbols fulfilling the following necessary conditions a \emph{candidate map}.

\begin{definition}\label{def:cand}
Let $G$ be the rotation group of a Platonic solid, and let $R(G)$ be its set of core rotations.  
Furthermore, let a list of orbit symbols of the two possible forms $(g)$ or $(g_1,g_2,g_3)$ with 
$g,g_1,g_2,g_3 \in  R(G)$, $g_3=g_1^{-1}g_2^{-1}$, be given such that the following conditions hold.
\begin{enumerate}
\item Each element $g\in R(G)$ appears at most once in an orbit symbol in the list.
\item (Circuit Property) The list of oriented combinatorial angles obtained from all the orbit symbols forms a single circuit of successively adjacent angles. (Two angles $(g_1,g_2)$, and $(g_3,g_4)$ are \emph{successively adjacent} if $g_3=g_2^{-1}$.) 
\item (Connectedness Property) $G$ is generated by the subset of $R(G)$ consisting of all group elements occurring in some orbit symbol in the list.
\end{enumerate}
Then the set of associated faces for each orbit symbol is constructed on a vertex set indexed by $G$ as follows, where $g^h=h^{-1}gh$.
An associated face of a face orbit $(g)$ of type $2$ is given by an oriented boundary walk through the ${\rm ord}(g)$ distinct vertices
 \[v_h, v_{hg^h}=v_{gh}, \ldots, v_{h(g^h)^{{\rm ord}(g)-1}}=v_{(g)^{{\rm ord}(g)-1}h}, \textnormal{ and back to } v_h \]
with dart labels $h^{-1}gh$, for each $h\in G$.
An associated face of a face orbit $(g_1,g_2,g_3)$ of type $1$ is given by an oriented boundary walk through the three distinct vertices
 \[v_h, v_{hg_2^h}=v_{g_2h}, v_{h(g_1^{-1})^h}=v_{g_1^{-1}h}, \textnormal{ and back to } v_h \]
with dart labels
\[h^{-1}g_2h,h^{-1}g_2^{-1}g_1^{-1}h,h^{-1}g_1h,\] 
 for each $h\in G$. Faces with the same boundary walks through the same vertices are identified. 
Edges result from identifying opposite darts with inverse labels between the same pair of vertices.
The polyhedral map resulting from the union of all associated faces of the orbits in the list is called a \emph{candidate map} for $G$.
\end{definition}

Evidently, all underlying maps of vertex-transitive polyhedra fulfill property 1. Then, the Circuit Property is equivalent to a disk neighborhood for each vertex (which ensures that we create a map on a closed surface), whereas the Connectedness Property is necessary and sufficient to ensure connectedness of the map. 

If the list of orbit symbols encodes the underlying map $M$ of a known vertex-transitive polyhedron, then we identify vertex $v_\mathbf{1}$ of the reconstructed candidate map with the initial vertex $v$ of the original map $M$.

\subsection{Geometric Isomorphism}\label{sec:background2:iso}

Assume that two maximally triangulated, vertex-transitive polyhedra $P, P'$ (vertex-transitive with respect to the same representation of the same group $G$) are congruent under an element $x$ in $O(3)$. Then, the underlying polyhedral maps $M$, $M'$ of $P$, $P'$ are isomorphic, and the labels of corresponding darts of the underlying \emph{labeled} maps are conjugate by $x$. The set of congruences yielding such a polyhedron $P'$ for $P$ forms a group, namely $N_{O(3)}(G)$, whose elements fix the system of axes for $G$ as a whole (preserving the order of each axis). Its normal subgroup $C_{O(3)}(G)$ fixes each axis individually. More specifically, we call the normal subgroup of index two consisting of the orientation-preserving symmetries of $\mathbb{E}^3$, \[N_{O(3)}^+(G)=N_{O(3)}(G)\cap SO(3)\cong N_{O(3)}(G)/C_{O(3)}(G)\] the \emph{geometric automorphism group} of $G$ because it is a subgroup of ${\rm Aut}(G)$ \cite{gallian}.

Given two isomorphic candidate maps $M$, $M'$ (by an isomorphism of polyhedral maps), we say that they are \emph{geometrically isomorphic} if the labels of corresponding darts are conjugate by the same element $x$ in $N_{O(3)}(G)$. Then, and only then, a realization of $M$ gives rise to a congruent realization of $M'$, where by realization we mean that $M$ is the underlying map of a polyhedron $P$. Geometric isomorphism of maps is an equivalence relation. Furthermore, if $x \in N_{O(3)}^+(G)$, then conjugation by $x$ of the labels of subsequent associated darts of a face $F$ of $M$ yield the labels of subsequent associated darts of a face $F'$ of $M'$. One could say that in this case the orientations of $M$ and $M'$ match up, whereas otherwise we need an additional change of orientation.

The notion of geometric isomorphism of maps induces the notion of geometric isomorphism of face orbits as given by their respective orbit symbol (which encodes both the dart labels and the orientation of faces).

\begin{definition}\label{isoclass}
Let $G$ be a Platonic rotation group, and let $R(G)$ be its set of core rotations. If $(g)$ is an orbit of type $1$ with ${\rm ord}(g)>2$, for $g\in R(G)$, then the \emph{geometric isomorphism class} of $(g)$ is defined as
\be\label{type1class} \left[(g)\right]_{N_{O(3)}(G)}\colon=\left\{(g)^A \vert A\in N_{O(3)}(G)\right\},\ee
where the orbit symbols $(g)^A$ are obtained by
\[(g)^A\colon=\begin{cases} (g^A)=(A^{-1}gA) &\textnormal{ for } A\in N_{O(3)}^+(G)\\
((g^{-1})^{A})=(A^{-1}g^{-1}A) &\textnormal{ for } A\in N_{O(3)}(G)\setminus N_{O(3)}^+(G).\end{cases} \]

If $(g,h,g^{-1}h^{-1})$ is an orbit of type $2$, for $g, h, g^{-1}h^{-1} \in R(G)$, then 
the \emph{geometric isomorphism class} of $(g,h,g^{-1}h^{-1})$ is defined as
\be\label{type2class}\left[(g,h,g^{-1}h^{-1})\right]_{N_{O(3)}(G)}\colon=\left\{(g,h,g^{-1}h^{-1})^A \vert A\in N_{O(3)}(G)\right\},\ee
where the orbit symbols $(g,h,g^{-1}h^{-1})^A$ are obtained by
\[(g,h,g^{-1}h^{-1})^A\colon=\begin{cases} (g^A, h^A, (g^{-1}h^{-1})^A) &\textnormal{ for } A\in N_{O(3)}^+(G)\\
((hg)^A,(h^{-1})^A,(g^{-1})^A) &\textnormal{ for } A\in N_{O(3)}(G)\setminus N_{O(3)}^+(G).\end{cases} \]
\end{definition}

We need only find one representative for each geometric isomorphism class of candidate maps and consider its realizability as a vertex-transitive polyhedron in order to classify vertex-transitive polyhedra of genus $\mathfrak{g}\geq 2$.

\section{Vertex-Transitive Polyhedra under Tetrahedral Rotation $T$}\label{sec:tet:intro}
\begin{figure}
\begin{center}
\includegraphics[width=.5\textwidth]{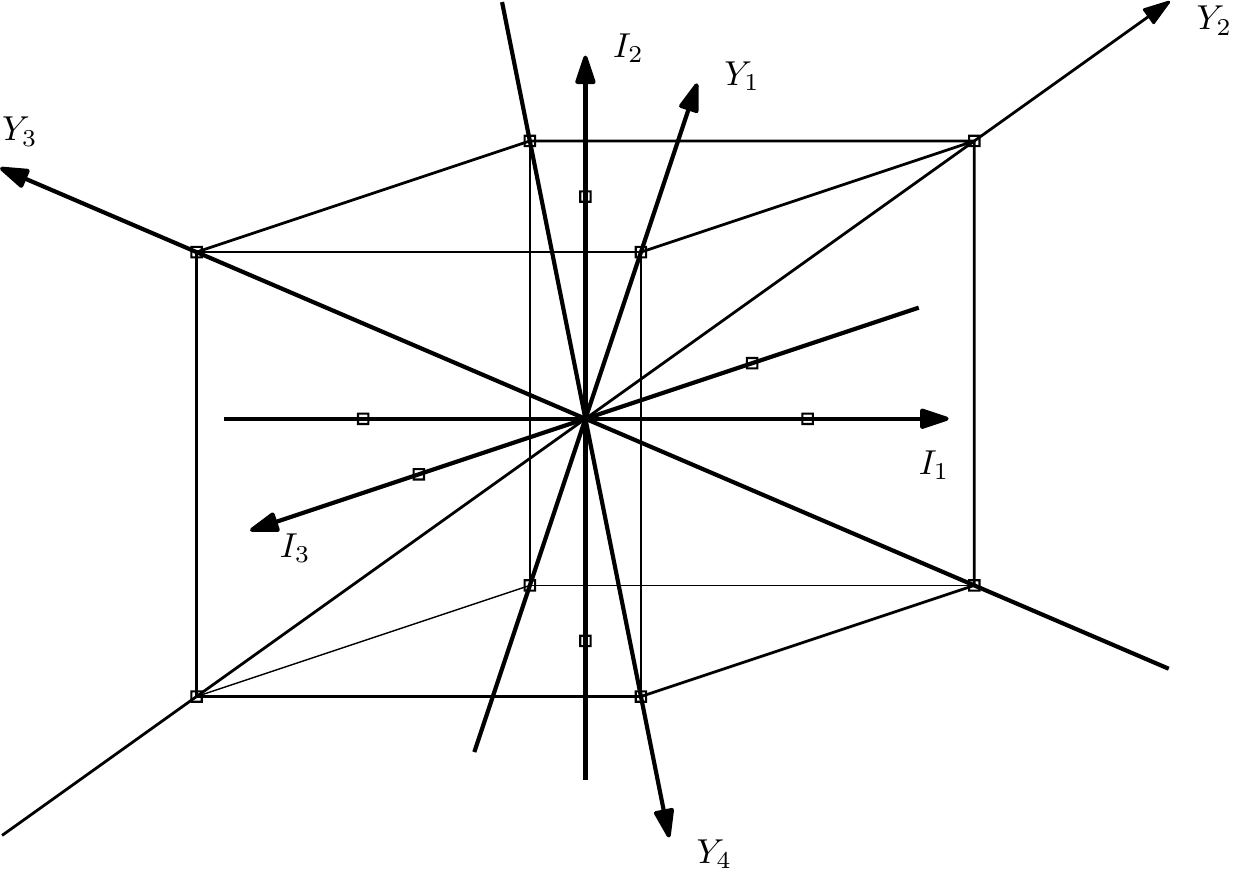}\\
\end{center}
\caption{The axes and corresponding rotation labels for the tetrahedral rotation group.}
\label{figtet}
\end{figure}

We begin our investigation by choosing, in Figure \ref{figtet}, the concrete configuration of axes in $\mathbb{E}^3$ for the tetrahedral rotation group $T$ of order $12$. 
Specifically, the directed axes corresponding to the involutory half-turns $I_1$, $I_2$, $I_3$ are also (in this order) the axes aligned with unit vectors $e_1=(1,0,0)$, $e_2=(0,1,0)$, $e_3=(0,0,1)$ of the coordinate axes. The axes of the three-fold rotations, denoted by the letter $Y$, point in the direction of the space diagonals of the auxiliary cube. 

The set of core rotations for $T$ consists of
\[R(T)=\{I_1,I_2,I_3,Y_1,Y_2,Y_3,Y_4,Y_1^{-1},Y_2^{-1},Y_3^{-1},Y_4^{-1}\}.\]
These group elements expressed as matrices are listed in Table \ref{tab:tetrot}. Multiplication is done from left to right. Note that for $T$ \emph{all} nontrivial elements are core rotations.
The geometric automorphism group $N^+_{O(3)}(T)$ of the tetrahedral rotation group $T$ is given by the octahedral group of rotations, and is isomorphic to $S_4$.

\begin{table}
\begin{tabular}{ccc}
$I_1=\left( \begin{array}{ccc}
1 & 0 & 0 \\
0 & -1 & 0 \\
0 & 0 & -1 \end{array} \right)$,
&
 $I_2=\left( \begin{array}{ccc}
-1 & 0 & 0 \\
0 & 1 & 0 \\
0 & 0 & -1 \end{array} \right)$,
&
$I_3=\left( \begin{array}{ccc}
-1 & 0 & 0 \\
0 & -1 & 0 \\
0 & 0 & 1 \end{array} \right)$,\\ 
&&\\
$Y_1=\left( \begin{array}{ccc}
0 & 1 & 0 \\
0& 0 & 1 \\
1& 0& 0 \end{array} \right)$,
&
$Y_2=\left( \begin{array}{ccc}
0 & 0 & -1 \\
1 & 0 & 0 \\
0 & -1 & 0 \end{array} \right)$,
&
$Y_3=\left( \begin{array}{ccc}
0 & 0 & -1 \\
-1 & 0 & 0 \\
0 & 1 & 0 \end{array} \right)$,\\
&&\\
$Y_4=\left( \begin{array}{ccc}
0 & 0 & 1 \\
-1 & 0 & 0 \\
0 & -1 & 0 \end{array} \right)$&&\\
\end{tabular}\\
\caption[The core rotations of $T$ as matrices.]{The core rotations of $T$ as matrices in the chosen representation (with respect to standard basis of $\mathbb{E}^3$.}
\label{tab:tetrot}
\end{table}

Maximal triangulation in the tetrahedral case leads necessarily to a proper triangulation of the surface or polyhedron. Thus, our focus lies on finding polyhedral realizations, free of non-trivial self-intersections, of triangulated surfaces with precisely $12$ vertices and tetrahedral symmetry. 

Edges occur in orbits under $T$ of size $6$ or $12$, depending on whether they are stabilized by an order~$2$ rotation or not. If edges intersect non-trivially, then, effectively, faces which contain these edges intersect. Thus we will focus on avoiding the non-trivial intersection of faces. 

As for faces, if a triangular face has a non-trivial stabilizer under the action of $T$ (face of type~$2$)
, then it is an equilateral triangle which is pierced in its center by an axis of rotation of order~$3$. The cyclic stabilizer of such a face has order~$3$ (generated by a rotation of order~$3$ around the piercing axis), and the face orbit will have size $\frac{12}{3}=4$. Recall that in this case each of the twelve vertices is incident to exactly one such face from the same face orbit. In the case of trivial face stabilizer under tetrahedral rotation symmetry (face of type~$1$), precisely three faces from the same face orbit are meeting at each vertex.

The genus of a surface with $12$ vertices is at most $6$, due to the Heawood inequality. 
Consequently, by Theorem \ref{thm:genus}, only genus $0$, $3$, $5$, $6$ can occur in a realizable map (a map embeddable as polyhedron).

\subsection{Face Orbits up to Geometric Isomorphism}\label{tetface}

Recall that simple vertex-transitivity allows us to derive an entire map $M$ from the local rotation at a chosen initial vertex $v$, or alternatively (and equivalently), from a list of orbit symbols (which we called candidate map).

Careful exhaustion of the possible selections of group elements in the tetrahedral case yields the following outcomes, up to geometric isomorphism for face orbits as defined in Definition \ref{isoclass}. 
First, the only possibility for (triangular!) faces of type 2 (with non-trivial stabilizer)
is given by the faces stabilized by a rotation of order $3$, with a corresponding orbit symbol $(g)$, where $g^3=\mathbf{1}$. Observe that all possibilities are conjugate by an element in $N_{O(3)}(T)$, which is the full octahedral group $\bar{O}$. 

If a face has trivial stabilizer, then the orbit symbol will be of the form $(g,h,g^{-1}h^{-1})$, with all entries being distinct core rotations. Now, if the orbit symbol does not contain an involution, there is only one possibility up to geometric isomorphism, of which $(Y_4^{-1}, Y_2^{-1}, Y_3^{-1})$ is representative. The reader may verify with the matrices in Table \ref{tab:tetrot} that $Y_4\cdot Y_2=Y_3^{-1}$.

If an orbit symbol contains at least two involutions, then the third element is also an involution, since $T$ contains the Klein Four group which is generated by each pair of the three involutions $I_1, I_2, I_3$. However, this kind of face orbit disqualifies itself, as it does not permit a connection to other orbits, and leads to a compound of (skew) tetrahedra (i.e., it can never be part of a candidate map). 

We conclude that the remaining possible orbit symbols for orbits of non-trivially stabilized faces contain exactly one involution and two elements of order $3$. It turns out that all these belong to the same class, $[(Y_1,Y_4,I_1)]_{N_{O(3)}(T)}$. 

In summary, there are precisely three geometric isomorphism classes of face orbits to consider: $[(Y_1)]_{N_{O(3)}(T)}$,  $[(Y_4^{-1}, Y_2^{-1}, Y_3^{-1})]_{N_{O(3)}(T)}$, and $[(Y_1,Y_4,I_1)]_{N_{O(3)}(T)}$.

\section{Enumeration of Candidate Maps}\label{sec:tet:maps}

The next step is to combine face orbits, encoded by orbits symbols, into a candidate map. 
Observe that the face orbits of type 2 (non-trivially stabilized) result simply from the missing inverses of elements in the list of orbits of type 1. Moreover, each map needs to contain at least one face orbit of type 1, and all face orbits of type 1 are connected to each other via matching inverse elements (successively adjacent combinatorial angles). We therefore focus on enumerating the possibilities for face orbits of type 1, and fill in the face orbits of type 2 as needed.

First, let us assume that the only orbits of type 1 (trivially stabilized faces) are of the same geometric isomorphism class as $(Y_4^{-1}, Y_2^{-1}, Y_3^{-1})$. W.l.o.g. then assume $(Y_4^{-1}, Y_2^{-1}, Y_3^{-1})$ is part of the selection of orbits. We obtain a first map by closing up the holes with orbits $(Y_2)$, $(Y_3)$, $(Y_4)$ of non-trivially stabilized faces. The genus is easily computed from the twelve vertices, the twelve faces in each orbit of type 1 and four in each orbit of type 2, and the number of edges as six times the number of elements in the list of orbits (each group element contributes a dart orbit, or half an edge orbit). Thus, our first completed map has genus $\mathfrak{g}=1-\frac{12-6\cdot 6+(12+3\cdot 4)}{2}=1$ and is therefore, by Theorem \ref{thm:genus}, not realizable as polyhedron (meaning it is not the underlying labeled, oriented map of some vertex-transitive polyhedron).

If other orbits of type 1 of the same class are involved, then precisely one such other orbit is involved (in order to avoid repetition of elements), and it must necessarily connect to $(Y_4^{-1}, Y_2^{-1}, Y_3^{-1})$ by precisely two matched up inverse elements. By geometric isomorphism, we may further require that the connection happens at $Y_2^{-1}$ and $Y_3^{-1}$, i.e. that the symbol of the new orbit contains the elements $Y_2$ and $Y_3$, which means that the new orbit is $(Y_2,Y_3,Y_1^{-1})$. However, matching inverse elements by adding in $(Y_4)$ and $(Y_1)$, we find that the circuit condition is violated; the subset of successively oriented combinatorial angles $(Y_4^{-1},Y_2^{-1})$, $(Y_2, Y_3)$, $(Y_3^{-1},Y_4^{-1})$, $(Y_4,Y_4)$ already gives a disk at $v$. Therefore, any other candidate maps which are geometrically non-isomorphic to our very first example of genus $1$ must contain a face orbit of type 1 from a different geometric isomorphism class also. 

Second, let us assume that the geometric isomorphism class of $(Y_1,Y_4,I_1)$ is allowed, and that there is at least one orbit from this class. W.l.o.g., we now assume the starting orbit for building our map to be $(Y_1,Y_4,I_1)$. We obtain a second candidate map, of genus $\mathfrak{g}=1-\frac{12-5\cdot 6+(12+2\cdot 4)}{2}=0$, by completing with $(Y_1^{-1})$ and $(Y_4^{-1})$. This is the map of the snub tetrahedron. 

Then, if no other orbit of type 1 from the same geometric isomorphism class is in the map, there is only the option of adding $(Y_4^{-1},Y_2^{-1},Y_3^{-1})$, and then the option of adding $(Y_3,Y_1^{-1}, Y_2)$ to that. In either case, the map created by filling in the missing orbits of type $2$ (none are missing for the latter!) has genus $4$ and is also excluded by Theorem \ref{thm:genus}. 

Next, still in the second case, suppose there is another orbit besides $(Y_1,Y_4,I_1)$ of the same class. Since any orbits geometrically isomorphic to $(Y_4^{-1}, Y_2^{-1}, Y_3^{-1})$ connect via inverses, so do the orbits geometrically isomorphic to $(Y_1,Y_4,I_1)$. We may therefore assume that the next added orbit of type 1 from this class connects to the element $Y_1$, i.e. that it contains the element $Y_1^{-1}$.  The options (recall that element $I_1$ cannot be repeated!) are \[ (Y_1^{-1}, Y_3^{-1}, I_2), (Y_1^{-1}, I_2, Y_2^{-1}), (Y_1^{-1}, Y_4^{-1}, I_3), (Y_1^{-1}, I_3, Y_3^{-1}).\] 
There are also eight more face orbits of the same class which contain neither $I_1$ nor $Y_1$ nor $Y_1^{-1}$:
\[(Y_2,I_3,Y_4^{-1}),(Y_2,Y_3^{-1},I_3),(Y_3,I_2,Y_4^{-1}), (Y_2,Y_4^{-1},I_2),\]
and the four orbits obtained by reversing the orientation, i.e. inverting each element and reversing the cyclic order of elements in the orbit symbol.

In this vein we continue to follow through the tree of possibilities, and leave it to the reader to verify that the lists of orbit symbols in Table \ref{tab:cand} give candidate maps, and all other possibilities violate one of the conditions (non-repetition of elements or Circuit Property) of Definition \ref{def:cand}.

\begin{table}
\begin{center}
\begin{tabular}{lllll}
\hline\noalign{\smallskip}
(a)&(b)&(c)&(d)&(e)\\
\noalign{\smallskip}\hline\noalign{\smallskip}
 $\begin{array}{l}
(Y_1 , Y_4 , I_1) \\
(Y_1^{-1} , I_2 , Y_2^{-1}) \\
(Y_4^{-1})\\
(Y_2) \end{array}$
&
 $\begin{array}{l}
(Y_1 , Y_4 , I_1) \\
(Y_1^{-1}, I_3, Y_3^{-1}) \\
(Y_4^{-1})\\
(Y_3) \end{array}$
&
 $\begin{array}{l}
(Y_1 , Y_4 , I_1) \\
(Y_1^{-1}, Y_3^{-1}, I_2) \\
(Y_4^{-1})\\
(Y_3) \end{array}$
&
 $\begin{array}{l}
(Y_1 , Y_4 , I_1) \\
(Y_1^{-1} , I_2 , Y_2^{-1}) \\
(Y_2 , Y_3^{-1} , I_3)\\
(Y_4^{-1})\\
(Y_3) \end{array}$\label{fig:cand:4}
 &
 $\begin{array}{l}
(Y_1 , Y_4 , I_1) \\
(Y_1^{-1}, Y_3^{-1}, I_2) \\
(Y_3, Y_2^{-1} , I_3)\\
(Y_4^{-1})\\
(Y_2) \end{array}$\\
\noalign{\smallskip}\hline
\end{tabular}
\end{center}
\caption{Candidate maps with orbits of class $[(Y_1,Y_4,I_1)]_{N_{O(3)}(T)}$ only.}
\label{tab:cand}
\end{table}

It is furthermore possible to show geometric isomorphism of the maps in Table \ref{tab:cand}(a) and \ref{tab:cand}(b) (by conjugation with the element $\mathbf{-1}\cdot Y_1^{-1}$ of $N_{O(3)}(T)$, where $\mathbf{-1}$ is the reflection in the origin), as well as geometric isomorphism of those in Table \ref{tab:cand}(d) and \ref{tab:cand}(e). We finish with four candidate maps which also fulfill the conditions on the genus, see Table \ref{tabletet}.

\begin{table}
\begin{center}
\begin{tabular}{p{1cm} p{0.5cm} p{2.5cm} p{7cm}}
\noalign{\smallskip}\hline
map & $\mathfrak{g}$	
 & type $1$ orbits & local rotation at $v$\\
 \noalign{\smallskip}\hline\noalign{\smallskip}
\begin{minipage}[b][1.5cm][c]{0.5cm}$M_0$\end{minipage}&\begin{minipage}[b][1.5cm][c]{0.5cm}$0$\end{minipage} 
&\begin{minipage}[b][1.5cm][c]{2.5cm}$(Y_1,Y_4,I_1)$\end{minipage}& \begin{minipage}[b][1.5cm][c]{7cm}$(Y_1, Y_1^{-1},Y_4,Y_4^{-1},I_1)$\end{minipage}\\

\begin{minipage}[b][2cm][c]{0.5cm}$M_1$\end{minipage}&\begin{minipage}[b][2cm][c]{0.5cm}$3$\end{minipage} 
&\begin{minipage}[b][2cm][c]{2.5cm} $\begin{array}{l}
(Y_1, Y_4 , I_1) \\
(Y_3^{-1}, Y_1^{-1}, I_3) \end{array}$\end{minipage} & \begin{minipage}[b][2cm][c]{7cm}$(Y_1,I_3,Y_3^{-1},Y_3,Y_1^{-1},Y_4, Y_4^{-1},I_1)$\end{minipage}\\

\begin{minipage}[b][2cm][c]{0.5cm}$M_2$\end{minipage}&\begin{minipage}[b][2cm][c]{0.5cm}$3$\end{minipage} 
&\begin{minipage}[b][2cm][c]{2.5cm} $\begin{array}{l}
(Y_1, Y_4,I_1) \\
 (Y_1^{-1}, Y_3^{-1} ,I_2) \end{array}$\end{minipage}&\begin{minipage}[b][2cm][c]{7cm}$(Y_1,Y_3^{-1},Y_3,I_2,Y_1^{-1},Y_4,Y_4^{-1},I_1)$\end{minipage}\\

\begin{minipage}[b][2cm][c]{0.5cm}$M_3$\end{minipage}&\begin{minipage}[b][2cm][c]{0.5cm}$6$\end{minipage} 
&\begin{minipage}[b][2cm][c]{2.5cm} $\begin{array}{l}
(Y_1 , Y_4 , I_1) \\
(Y_1^{-1} , I_2 , Y_2^{-1}) \\
(Y_2 , Y_3^{-1} , I_3) \end{array}$ \end{minipage}&\begin{minipage}[b][2cm][c]{7cm}$(Y_1,I_2,Y_2^{-1},Y_3^{-1},Y_3,I_3,Y_2,Y_1^{-1},Y_4,Y_4^{-1},I_1)$\end{minipage}\\
\noalign{\smallskip}\hline
\end{tabular}
\end{center}
\caption{Remaining candidate maps for the tetrahedral case, up to geometric isomorphism.}
\label{tabletet}
\end{table}

We see that $\mathfrak{g}=0$ is not of higher genus; it is the map of the snub tetrahedron. The case $\mathfrak{g}=6$ for $12$ vertices has been ruled out in full generality (for any polyhedral map, without symmetry assumptions) by the work of Schewe \cite{schewe}. Hence, map $M_3$ is also irrelevant for our search for vertex-transitive polyhedra of higher genus with tetrahedral symmetry.

\begin{figure}
\begin{center}
\includegraphics[width=.4\textwidth]{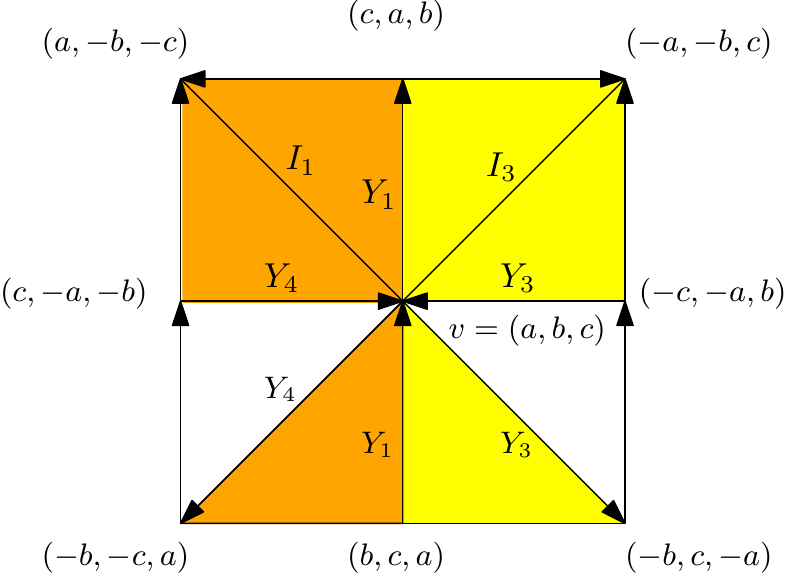}
\end{center}
\caption{The vertex-star of map $M_1$ of genus $3$.}
\label{map1}
\end{figure}

\begin{figure}
\begin{center}
\includegraphics[width=.4\textwidth]{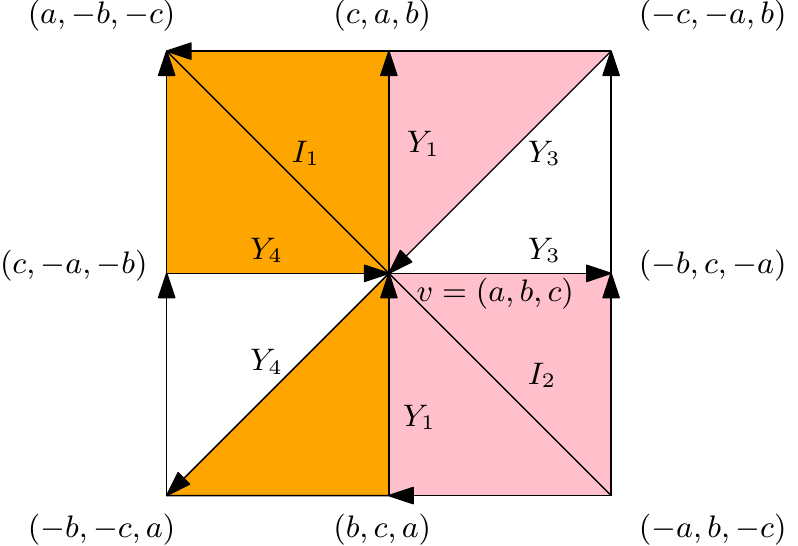}
\end{center}
\caption{The vertex-star of map $M_2$ of genus $3$.}
\label{map2}
\end{figure}

$M_1$ and $M_2$ of genus $3$, are of Schl\"afli type $\{3,8\}$, with $8$ triangles meeting at each vertex. Their initial vertex stars, together with vertex coordinates, are depicted in Figures \ref{map1} and \ref{map2}, where the coordinates of the initial vertex $v=v_{\mathbf{1}}$ are assumed to be $(a,b,c)$. These two maps are related, as a combinatorial operation known as \emph{edge flip} performed on the first map $M_1$ yields the second map $M_2$. We will show that only $M_1$ is realizable as a polyhedron in $\mathbb{E}^3$.

\section{Reduction to a Single Feasible Map}\label{sec:tet:reduction}
A candidate map $M$ has a vertex set indexed by $T$, with initial vertex $v_{\mathbf{1}}$, at which the orbit symbols for the face orbits can be read off. Using the abbreviation $v=v_{\mathbf{1}}$, the question becomes whether there is an assignment of coordinates $v=(a,b,c)$ which realizes $M$ as a vertex-transitive polyhedron.
Certainly, $a=b=c=0$ does not lead to a polyhedron. Then, w.l.o.g., $v$ lies on the unit sphere around the origin, so $a^2+b^2+c^2=1$ (otherwise $P$ can be stretched or shrunk accordingly).

For map $M_2$, pictured in Figure \ref{map2}, we now prove that such a realization as a vertex-transitive polyhedron does not exist. The following obstruction plays a key role. 
\emph{When an axis of rotation pierces a triangle or an edge which its rotations do not stabilize, then a self-intersection is produced within the orbit of the face under the tetrahedral rotation group.} In this case we do not have a polyhedron. In terms of determinants, if $p$ is a row vector aligned with the axis in question, and $v_0$, $v_1$, $v_2$ are row vectors giving the vertices of a triangle face in the map, then the axis pierces triangle $v_0v_1v_2$ in its relative interior or on its relative boundary if and only if no two determinants out of the three determinants
\[|p^T\ v_0^T\ v_1^T|, |p^T\ v_1^T\ v_2^T|, |p^T\ v_2^T\ v_0^T|\]
have opposite signs (i.e., one is positive and one is negative).

We will show that one of the following cases occur for any assignment of vertex coordinates for map $M_2$:
\begin{paragraph}{Case 1} The triangle $\Delta_1$, stabilized by $\langle Y_3 \rangle$, with vertices $v=(a,b,c)$, $(-b,c,-a)$, and $(-c,-a,b)$ is pierced by the axis of $I_3$ (and also by the symmetric axes of $I_1$ and $I_2$); this happens precisely when the three determinants
\[\left|\begin{array}{ccc}0 & a & -b\\
0 & b & c\\
1 & c &-a
\end{array}\right|=b^2+ac,\]
\[\left|\begin{array}{ccc}0 & -b & -c\\
0 & c & -a\\
1 & -a &b
\end{array}\right|=c^2+ab,\]
and 
\[\left|\begin{array}{ccc}0 & -c & a\\
0 & -a & b\\
1 & b & c
\end{array}\right|=a^2-bc\] 
are all non-negative or all non-positive. Observe that if any two of these determinants are negative, then the third must be positive, 
thus we may as well say that this case happens only when all three determinants are non-negative.
\end{paragraph}
\begin{paragraph}{Case 2} The triangle $\Delta_2$, stabilized by $\langle Y_4 \rangle$, with vertices $v=(a,b,c)$, $(c,-a,-b)$, and $(-b,-c,a)$ is pierced by the axis of $I_3$ (and also by the symmetric axes of $I_1$ and $I_2$); this happens precisely when the three determinants
\[\left|\begin{array}{ccc}0 & a & c\\
0 & b & -a\\
1 & c &-b
\end{array}\right|=-(a^2+bc),\]
\[\left|\begin{array}{ccc}0 & c & -b\\
0 & -a & -c\\
1 & -b & a
\end{array}\right|=-(c^2+ab),\]
and 
\[\left|\begin{array}{ccc}0 &-b & a\\
0 & -c & b\\
1 & a &c
\end{array}\right|=-(b^2-ac)\]
are all non-negative or all non-positive. Observe that if any two of these determinants are positive, then the third must be negative, 
thus we may as well say that this case happens only when all three determinants are non-positive.
\end{paragraph}
\begin{paragraph}{Case 3} The triangle $\Delta_3$, of the orbit $(Y_1, Y_4, I_1)$, with vertices $v=(a,b,c)$, $(-b,-c,a)$, and $(b,c,a)$ is pierced by the axis of $Y_4$; 
this happens precisely when the three determinants
\[\left|\begin{array}{ccc}1 & a & -b\\
-1 & b & -c\\
1 & c &a
\end{array}\right|=\frac{(a+b)^2}{2}+\frac{(b+c)^2}{2}+\frac{(c-a)^2}{2}\geq 0,\]
\[\left|\begin{array}{ccc}1 & -b & b\\
-1 & -c & c\\
1 & a &a
\end{array}\right|=-2a(b+c),\]
and 
\[\left|\begin{array}{ccc}1 & b & a\\
-1 & c & b\\
1 & a &c
\end{array}\right|=-(b+c)(a-c)+(b+a)(b-a)=\underbrace{(b+c)^2}_{\geq 0}-(a+c)(b+a)\]
are all non-negative.\end{paragraph}
\begin{paragraph}{Case 4} The triangle $\Delta_4$, of the orbit $(Y_1^{-1}, Y_3^{-1}, I_2)$, with vertices $v=(a,b,c)$, $(-c,-a,b)$, and $(c,a,b)$ is pierced by the axis of $Y_3$; this happens precisely when the three determinants
\[\left|\begin{array}{ccc}-1 & a & -c\\
1 & b & -a\\
1 & c & b
\end{array}\right|=-\frac{(a+b)^2}{2}-\frac{(c-b)^2}{2}-\frac{(c+a)^2}{2}\leq 0,\]
\[\left|\begin{array}{ccc}-1 & -c & c\\
1 & -a & a\\
1 & b &b
\end{array}\right|=2b(a+c),\]
and 
\[\left|\begin{array}{ccc}-1 & c & a\\
1 & a & b\\
1 & b &c
\end{array}\right|=(b-c)(a+c)-(a+b)(a-b)=\underbrace{-(a+c)^2}_{\leq 0}+(a+b)(b+c)\]
are all non-positive.
\end{paragraph}
\begin{paragraph}{Case 5} The triangles $\Delta_3$ and $\Delta_4$ intersect.
\end{paragraph}

\medskip
Now, the search space for $(a,b,c)$ splits into regions which can be matched up with the above cases as follows.
\begin{paragraph}{Region 1} We have $c^2+ab<0$. This means that $a, b$ must have different signs and cannot be zero, and either $|a|\geq |b|,|c|$ or $|b|\geq |a|,|c|$. Now, if $|a|\geq |b|, |c|$ then Case 4 occurs, as $a\pm b$, $a \pm c$ have the same sign as $a$ or are zero, while one of $b \pm c$ must have the same sign as $b$. Otherwise, if $|b|\geq |a|, |c|$ then Case 3 occurs, as $b\pm a$, $b \pm c$ have the same sign as $b$ or are zero, while one of $a \pm c$ must have the same sign as $a$.\end{paragraph}
\begin{paragraph}{Region 2} We have $a^2+bc < 0$ and $b^2+ac < 0$. In this case, $a$, $b$ are not zero, have a sign different from $c\neq 0$, and $|c|\geq |a|,|b|$. It is easily checked that this satisfies the conditions for both Case 3 and Case 4.
\end{paragraph}
\begin{paragraph}{Region 3} We have $a^2-bc < 0$ and $b^2-ac < 0$. It follows that $a,b,c$ are not zero, have the same sign, and $|c|> |a|, |b|$. We now show that the triangle $\Delta_3$ of $M_2$ with vertices $v=(a,b,c)$, $(-b,-c,a)$, $(b,c,a)$ (see Figure~\ref{map2}) intersects the triangle $\Delta_4$ with vertices $v=(a,b,c)$, $(-c,-a,b)$, and $(c,a,b)$ in a relative interior point (satisfying Case 5 above). Note that each of these triangles has an edge (the one opposite to $v$) which is fixed by $I_3$, and therefore pierced perpendicularly by the corresponding axis in its center. It is these edges that will be of interest. Specifically, the determinant 

\begin{eqnarray*}\left\vert\begin{array}{cccc}
1&1&1&1\\
-c&a&-b&b\\
-a&b&-c&c\\
b&c&a&a
\end{array}\right\vert&=&-2 (c-b) (a^2+a b-b^2+a c-b c-c^2)\\
&=&-2(c-b)\underbrace{\left(2(a^2-bc)-\frac{(a-b)^2}{2}-\frac{(c-a)^2}{2}-\frac{(c-b)^2}{2}\right)}_{<0}
\end{eqnarray*}

has the same sign as $c$, and the determinant 

\begin{eqnarray*}
\left\vert\begin{array}{cccc}
1&1&1&1\\
c&a&-b&b\\
a&b&-c&c\\
b&c&a&a
\end{array}\right\vert&=&-2 (c+b) (a^2-a b+b^2-a c-b c+c^2)\\
&=&-2(c+b)\underbrace{\left(\frac{(a-b)^2}{2}+\frac{(c-a)^2}{2}+\frac{(c-b)^2}{2}\right)}_{> 0}
\end{eqnarray*}

has the opposite sign as $c$. 

\pagebreak
Similarly, the determinant

\begin{eqnarray*}
\left\vert\begin{array}{cccc}
1&1&1&1\\
-b&a&-c&c\\
-c&b&-a&a\\
a&c&b&b
\end{array}\right\vert&=&2 (c-a) (-a^2+a b+b^2-a c+b c-c^2)\\
&=&2(c-a)\underbrace{\left(-2(-b^2+ac)-\frac{(a-b)^2}{2}-\frac{(c-a)^2}{2}-\frac{(c-b)^2}{2}\right)}_{< 0}
\end{eqnarray*}

has the opposite sign as $c$, and the determinant
\begin{eqnarray*}
\left\vert\begin{array}{cccc}
1&1&1&1\\
b&a&-c&c\\
c&b&-a&a\\
a&c&b&b
\end{array}\right\vert&=&2 (a+c) (a^2-a b+b^2-a c-b c+c^2)\\
&=&2(c+a)\underbrace{\left(\frac{(a-b)^2}{2}+\frac{(c-a)^2}{2}+\frac{(c-b)^2}{2}\right)}_{> 0}
\end{eqnarray*}
has the same sign as $c$. This means that there is a non-trivial intersection between the two triangles $\Delta_3$ and $\Delta_4$, as the plane of either triangle strictly separates the vertices of the edge opposite $v$ in the other triangle. (Both triangles must lie in the same closed half-space, since all vertices lie on a sphere.)
\end{paragraph}
\begin{paragraph}{Region 4} In the remaining region, which is disjoint from the previous ones, the inequality $c^2+ab\geq 0$ holds. Furthermore, at least one expression in each of the following four pairs is nonnegative: $a^2+bc$ and $b^2+ac$,  $a^2-bc$ and $b^2-ac$, $b^2\pm ac$, $a^2\pm bc$. Therefore, the conditions for Case 1 or Case 2 above are fulfilled. 
\end{paragraph}

Consequently, all coordinate vectors are infeasible and the second map $M_2$ of genus $3$ is not realizable as a polyhedron.

\section{The Unique Vertex-Transitive Polyhedron of Genus $3$ under $T$}\label{sec:tet:conclusion}
\begin{table}[ht]
\centering
\subtable[Vertex coordinates.]{
\begin{tabular}{c|rrr}
vertex& $x_1$&$x_2$&$x_3$\\\hline
$v=v_{\mathbf{1}}$&$1$&$ 2$&$ 6$\\
$v_2$&$-2$&$6$&$-1$\\
$v_3$& $-6$&$-1$&$2$\\
$v_4$&$6$&$1$&$ 2$\\
$v_5$&$-1$& $-2$& $6$\\
$v_6$&$2$&$-6$& $-1$\\
$v_7$&$1$ & $-2$ & $-6$\\
$v_8$&$-2$&$ -6$&$1$\\
$v_9$&$-6$&$1$&$ -2$\\
$v_{10}$&$-1$&$2$&$-6$\\
$v_{11}$&$2$&$ 6$&$ 1$\\
$v_{12}$&$6$&$ -1$&$ -2$\\
\end{tabular}}\quad\quad
\subtable[List of triangles for $M_1$.]{
\begin{tabular}{lll}
$v_1v_2v_3$ &$v_1v_3v_5$ &$v_1v_4v_7$\\
$v_4v_5v_6$&$v_1v_5v_4$&$v_1v_7v_{12}$\\
$v_7v_8v_9$&$v_2v_1v_{11}$&$v_2v_{10}v_{6}$\\
$v_{10}v_{11}v_{12}$&$v_2v_{11}v_{10}$&$v_2v_{6}v_7$\\
$v_1v_{12}v_8$&$v_3v_2v_{9}$&$v_3v_8v_{12}$\\
$v_2v_7v_{4}$&$v_3v_{9}v_{8}$&$v_3v_{12}v_6$\\
$v_3v_6v_{10}$&$v_4v_6v_{12}$&$v_4v_{11}v_{9}$\\
$v_5v_{9}v_{11}$&$v_4v_{12}v_{11}$&$v_4v_{9}v_2$\\
& $v_6v_5v_{8}$&$v_5v_3v_{10}$\\
&$v_6v_8v_{7}$&$v_5v_{10}v_9$\\
&$v_7v_9v_{10}$&$v_8v_5v_{11}$\\
&$v_7v_{10}v_{12}$&$v_8v_{11}v_{1}$\\					
\end{tabular}
}
\caption{Data for the polyhedron of genus $\mathfrak{g}=3$ with underlying map $M_1$.}
\label{tab:map1coord}
\end{table}

\begin{figure}[h]
\begin{center}\includegraphics[width=0.5\textwidth]{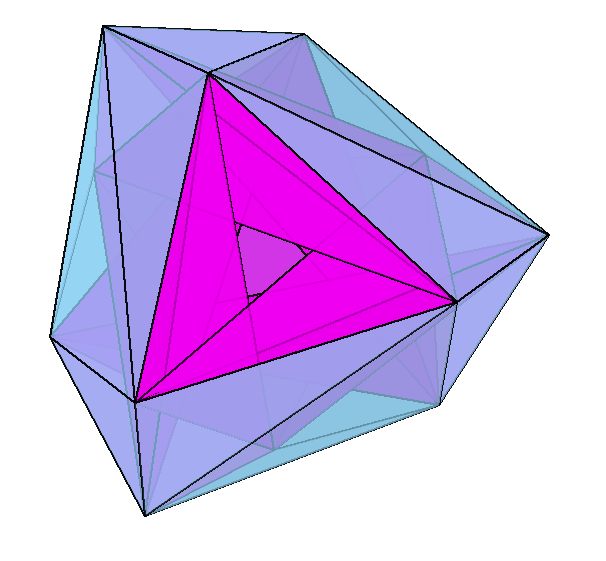}\end{center}
\caption[The only higher genus polyhedron with vertex-transitivity under $T$.]{The only higher genus polyhedron with vertex-transitivity under the tetrahedral group of rotations.}
\label{fig:tet_poly}
\end{figure}

The existence of a vertex-transitive polyhedron of genus $3$ and type $\{3,8\}$ under tetrahedral symmetry has already been postulated by Gr\"unbaum and Shephard \cite{gs1984} but no coordinates were given. A construction was also given in \cite{sw1986}. Here we are able to supply integer coordinates for the base vertex of the remaining candidate map $M_1$, namely $v=(1,2,6)$.
The local rotation, i.e. the sequence of labels of outgoing darts in counterclockwise orientation around $v$ is $(Y_1,I_3,Y_3^{-1},Y_3,Y_1^{-1},Y_4, Y_4^{-1},I_1)$ (compare also Figure \ref{map1}). Recall that the labels stand for the smallest non-trivial rotations (core rotations $R(T)$) around the axes of the tetrahedral group $T$ as represented in Figure \ref{figtet}, where the auxiliary cube is aligned with the coordinate axes of $\mathbb{E}^3$. Refer to Table \ref{tab:map1coord} for a complete list of vertex coordinates for the $12$ vertices and a list of the $32$ triangle faces.

\begin{figure}[h]
\begin{center}\includegraphics[width=0.45\textwidth]{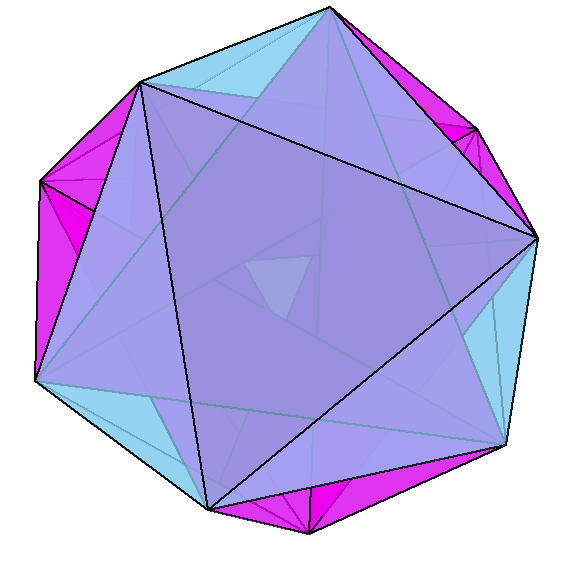}\end{center}
\caption[The only higher genus polyhedron with vertex-transitivity under $T$.]{The only higher genus polyhedron with vertex-transitivity under the tetrahedral group of rotations.}
\label{fig:tet_poly2}
\end{figure}

In Figures \ref{fig:tet_poly} and \ref{fig:tet_poly2}, we see two images of the polyhedron as constructed with the software JAVAVIEW~\cite{polthier:javaview}. Its convex hull is combinatorially isomorphic to a snub tetrahedron. It is impossible to realize this triangulated polyhedron with some coplanar faces, as we will see now by inspecting the map~$M_1$ (see Figure~\ref{map1}). 

For a non-maximally triangulated version, we would have to combine several faces into a larger face by making them coplanar, while also respecting the symmetry. However, for map $M_1$, it is impossible to incorporate the faces adjacent to a given non-trivially stabilized face, such as the ones symmetrically surrounding the face at $v$ stabilized by $Y_1$; some diagonals of this putative hexagon are already part of the polyhedron and would therefore cause a self-intersection. It is also impossible to make two adjacent, trivially stabilized faces of the same orbit coplanar. If this were successful, we could re-create a maximal triangulation, this time using the other possible splitting edge, which would yield a realization of the (infeasible!) map $M_2$ of genus~$3$. Clearly, this is a contradiction. The last possibility is coplanarity of adjacent, trivially stabilized faces from different orbits. We can exclude this possibility as well, since a convex face thus created cannot have two edges which are pierced by order $2$ axes, without the face being pierced by the third order $2$ axis in its interior. 

We thus conclude the proof of Theorem \ref{thm:tet}; there is only combinatorial type of polyhedron with vertex-transitive tetrahedral symmetry and positive genus. We observe that a left-handed and a right-handed version is possible. Note that neither map $M_1$ (of the realizable polyhedron) nor map $M_2$ are combinatorially isomorphic to Dyck's regular map (compare \cite{sw1986}, \cite{brehm1987}).

\section{Notes}\label{sec:closing}
The combinatorial types of vertex-transitive polyhedra of genus $\mathfrak{g}\geq2$ under tetrahedral symmetry have been exhausted, with a slightly different approach as presented in the author's dissertation \cite{leopold2014}. The remaining cases of octahedral and icosahedral symmetry are much more involved. Progress and partial results with a focus on gaining geometric insight into the problem will be published in another paper \cite{leopold2015}. 

We have focused on the genus range $\mathfrak{g}\geq 2$ in this article. However, note that for toroidal vertex-transitive polyhedra, we must necessarily have simple transitivity of the symmetry group as well. This follows since reflections are also not allowed for $\mathfrak{g}=1$ (see \cite{gsw2013}), the torus cannot have Platonic symmetry at all (Theorem \ref{thm:genus}), and since a polyhedron cannot live in a plane. Infinite two-parameter families for the remaining possible symmetry groups of dihedral type have already been constructed, first in \cite{gs1984} and, independently, by Brehm (private communication), and described in more detail in \cite{gsw2013}. It appears unlikely that more variations could be found.

We conclude this paper with two brief remarks on ideas which emerged in the process of structuring the problem. 
First, there is a nice way to determine face orbits geometrically by \emph{pole figures}, figures (single points or spherical triangles) created by the intersection points of the positive axis directions for the group elements in the orbit symbol and the unit sphere. Specifically, for orbit $(g_1,g_2,g_3)$, the corresponding triangle of poles can be seen as a spherical triangle with interior angle $\pi-\pi/{\rm ord}(g)$ at the pole for $g$. Not only is this idea helpful in classifying possible face orbits, it also visualizes the symmetry of each candidate map, even when we do not yet know whether it is realizable as a vertex-transitive polyhedron. Geometric isomorphism, as defined in Section \ref{sec:background2:iso} simply means congruence (under an element in $N_{O(3)}(T)$) for pole figures. 

Second, we observe that the circuit condition of Definition \ref{def:cand} can also be verified diagrammatically. In order to do this, depict each orbit of type 1 as a triangle with the entries of the orbit symbol labeling the oriented boundary walk, and an orbit of type 2 as a monogon whose oriented boundary walk (just one dart) is again labeled with the entry. Then the circuit condition checks whether gluing of these faces at opposite darts, respectively identifying head and tail vertex of any dart labeled with an involution, also identifies all vertices in the resulting closed, oriented surface. It is this interpretation which makes the connection to branched covers, quotient surfaces, and embedded voltage graphs.

\section*{Acknowledgements}
Egon Schulte introduced me to this interesting and challenging problem when I was his Ph.D. student at Northeastern University. I am very grateful for his guidance, encouragement and many useful discussions.
\bibliographystyle{plain}
\bibliography{bibliography3}

\begin{flushleft}
Undine Leopold\\
Technische Universit\"at Chemnitz \\
Fakult\"at f\"ur Mathematik\\
D - 09107 Chemnitz\\
Germany\\[0.5ex]
undine.leopold@mathematik.tu-chemnitz.de
\end{flushleft}
\end{document}